\documentclass[12pt]{amsart}
\renewcommand{\Bbb}{\mathbb}

\newcommand{\Zzz}{\Bbb Z}

\newtheorem{theorem}{Theorem}

\newtheorem{lemma}{Lemma}

\newtheorem{definition}{Definition}
\newtheorem{corollary}{Corollary}

\newcommand{\vanish}[1]{}

\begin{document}

\title
[A Combinatorial proof of a result of Hetyei and Reiner]
{A Combinatorial proof of a result of Hetyei and Reiner on 
Foata-Strehl type permutation trees}
            
\author{Mikl\'os B\'ona}
\thanks{This paper was written while the author was
a one-term visitor at  Mathematical Sciences Research Institute in Spring 
1997. This visit
was supported by an MIT Applied Mathematics Fellowship.}
\address{Department of Mathematics \\
        Massachusetts Institute of Technology \\
        Cambridge, MA 02139}
            
\begin{abstract} We give a combinatorial proof of the known \cite{hetyei}
 result that
there are exactly $n!/3$ permutations of length $n$ in the minmax tree 
representation of which the $i$th node is a leaf. We also prove the
new result that the number of $n$-permutations in which this node has one
child is $n!/3$ as well, implying that the same holds for those in which
this node has two children. \end{abstract}

\maketitle  

\section{Background and definitions}
In \cite{hetyei} the authors gave a new proof for a theorem of Purtill 
\cite{purtill} by introducing a new $\Zzz_2^{n-1}$-action on the 
symmetric group
$S_n$. This had the flavor of a $\Zzz_2^n$-action defined and succesfully used
by Foata and Strehl \cite{Foata-Strehl1} \cite{Foata-Strehl2},
 but it proved to have a more symmetric structure.
 In order to define their new group action, Hetyei and Reiner studied a 
binary tree representation of the
 permutation $p=(p_1,p_2,\cdots ,p_n)$ called  {\em minimum-maximum trees}
or, in the rest of this paper, {\em minmax} trees. These are defined as
 follows.
\begin{definition} Let $p$ be a permutation of length $n$ and
let $p = u\,m\,v$ where $m$ is the leftmost of the
 minimum and maximum
letters of $p$,
$u$ is the subword preceding $m$  and $v$ is the subword  following 
$m$.
The {\em minmax tree} $T^m_p$ has $m$ as its
root.  The right subtree of $T^m_p$ is obtained by applying the
definition recursively to $v$. 
Similarly, the left subtree of $T^m_p$ is obtained by applying the
definition recursively to $u$. \end{definition}
So every node of $T^m_p$ which is not a leaf is either a minimum node or a
 maximum node. Useful information can be read off these trees, for example,
it can be shown that $p$ is an Andr\'e-permutation if and only
if all these nodes are minimum nodes.

 The minmax tree of the permutation  $p = 3\,6\,7\,1\,5\,2\,10\,4\,9\,8$ 
is shown on Figure 1.

$$
 \begin{picture}(378,122)(105,677)
\thicklines
\put(160,760){\circle*{6}}
\put(240,780){\circle*{6}}
\put(320,760){\circle*{6}}
\put(280,740){\circle*{6}}
\put(200,740){\circle*{6}}
\put(240,720){\circle*{6}}
\put(360,740){\circle*{6}}
\put(400,720){\circle*{6}}
\put(440,700){\circle*{6}}
\put(480,680){\circle*{6}}
\put(160,760){\line( 4, 1){ 80}}
\put(240,780){\line( 4,-1){ 80}}
\put(320,760){\line( 2,-1){160}}
\put(320,760){\line(-2,-1){ 40}}
\put(160,760){\line( 2,-1){ 80}}
\put(155,770){3}
\put(235,790){1}
\put(315,770){2}
\put(275,750){5}
\put(195,750){6}
\put(235,730){7}
\put(355,750){10}
\put(395,730){4}
\put(435,710){9}
\put(475,690){8} 
\put(175,660){Figure 1: $T^m_p$ for  $p = 3\,6\,7\,1\,5\,2\,10\,4\,9\,8$.}
\end{picture}
$$

\vskip .8 cm

If $a=a_1,a_2,\cdots ,a_k$ is a subword of the $n$-permutation $p$, then the
{\em pattern} of $a$ is the unique $k$-permutation for which the pairwise
comparisons between the $i$th and $j$th entries are the same as those for
$a$. For example, the pattern of the subword  $3\,6\,7\,1$ is $2\,3\,4\,1$,
while the paattern of the subword $4\,9\,8$ is $1\,3\,2$.

The group action defined by Hetyei and Reiner is as follows: Let $1\leq i
\leq n$. Then $\psi_i(T^m_p)$ is obtained from $T^m_p$ by changing only the
subtree whose root is the entry $p_i$. If $p_i$ was a minimum node, then we
take the {\em maximal} entry among all entries of this subtree, put it into
the place of $p_i$ (so to the root of this subtree), and write the other
entries of the subtree to the nodes of the subtree so that their {\em pattern}
is the same as it was before the operation. Similarly, if $p_i$ was a maximum
node, then we put the {\em minimal} among the entries of its subtree to the
place of $p_i$, then we write the other entries in the subtree to the nodes
so that their pattern is unchanged. 

Figure 2 shows the image of our previuos example 
 $p$, under the operator $\psi_7$.

$$
 \begin{picture}(378,122)(105,677)
\thicklines
\put(160,760){\circle*{6}}
\put(240,780){\circle*{6}}
\put(320,760){\circle*{6}}
\put(280,740){\circle*{6}}
\put(200,740){\circle*{6}}
\put(240,720){\circle*{6}}
\put(360,740){\circle*{6}}
\put(400,720){\circle*{6}}
\put(440,700){\circle*{6}}
\put(480,680){\circle*{6}}
\put(160,760){\line( 4, 1){ 80}}
\put(240,780){\line( 4,-1){ 80}}
\put(320,760){\line( 2,-1){160}}
\put(320,760){\line(-2,-1){ 40}}
\put(160,760){\line( 2,-1){ 80}}
\put(155,770){3}
\put(235,790){1}
\put(315,770){2}
\put(275,750){5}
\put(195,750){6}
\put(235,730){7}
\put(355,750){4}
\put(395,730){8}
\put(435,710){10}
\put(475,690){9} 
\put(165,660){Figure 2: $\psi_7(T^m_p)$ for 
 $p = 3\,6\,7\,1\,5\,2\,10\,4\,9\,8$.}
\end{picture}
$$

\vskip .8 cm

It is clear that $\psi_i(T^m_p)=T^m_p$ if and only if $p_i$ is a leaf in 
$T^m_p$. 
Hetyei and Reiner  proved that any for any $i\leq n-2$, $p_i$ is a leaf
 in exactly
$n!/3$ permutations of length $n$ (in the rest of this paper, 
$n$-permutations), thus $\psi_i$ fixes $n!/3$ $n$-permutations.
 This is a little surprising: it is not even intuitively 
obvius why this number should be independent from $i$, and indeed, in the
original group action of Foata and Strehl, it was not. Hetyei and Reiner used 
exponential generating functions 
and partial 
differential equations to get this result, however they noted that such a
 simple formula would
deserve a combinatorial proof. 

Present paper provides such a proof. We also
prove the new result that the number of permutations in which $p_i$ has 1 child
as well as the number of $n$-permutations in which $p_i$ has two children is
$n!/3$, too, whenever $2\leq i \leq n-2$. Thus this group action has even more
symmetries, which supports the inventors' claim that it could be used in
analysing random permutations.

\section{The Proof of the Theorem}

First we make some simple observations which will be useful later.
For any $i\leq n-1$, either $p_i$ is an ancestor of $p_{i+1}$, or $p_{i+1}$ is
 an  ancestor of $p_i$, for otherwise these two entries would have a common
 ancestor,
which would have an index {\em between} $i$ and $i+1$, a contradiction. So in
particular, no permutation's minmax tree can have both $p_i$ and $p_{i+1}$ as
leaves. This allows the following definition.

\begin{definition} Let $1\leq i\leq n-2$. Then the {\em $i$th local extremum}
of a permutation $p$ is the entry which is closest to the root of the minmax
tree of $p$ among $p_i$, $p_{i+1}$ and $p_{i+2}$. This entry will be denoted
$e_i$. \end{definition}

Being closest to the root means having the shortest path to the root in this 
definition.
Note that $e_i$ always exists: either $p_{i+1}$ is an ancestor of both of its
neighbors and $e_i=p_{i+1}$, or only one of its neighbors, say $p_i$ is an 
ancestor of $p_{i+1}$ and $p_{i+2}$ is a descendant of $p_{i+1}$, implying
$e_i=p_i$, or both of its neighbors are ancestors of $p_{i+1}$, in which case,
by the tree-property, one must be an ancestor of the other.

\begin{lemma} The number of $n$-permutations whose minmax tree contains $p_1$
as a leaf is $n!/3$. \end{lemma}
\begin{proof}
By symmetry, we can suppose that the entry 1 of our permutation precedes the 
entry $n$. Indeed, the complement of $p$ (the permutation whose $i$th entry is
$n+1-p_i$) has an isomorphic minmax tree to $T^m_p$.

 First we consider the case when
the entry 1  is among the first three elements, so in
particular, $e_1=1$. This gives rise to three subcases:
\begin{itemize}
\item If $p_1=1$, then $p_1$ is the root of the minmax tree.
\item If $p_2=1$, then $p_2$ is the root, and $p_1$ is a leaf.
\item If $p_3=1$, then $p_3$ is the root, its left subtree has
$p_1$ and $p_2$ as nodes, and among these, by definition, $p_2$ is a leaf, and
 $p_1$ is not.
\end{itemize}

Clearly, these subcases are equally likely to occur, so each of them 
occurs with probability 1/3.  

Now suppose that the entry of the permutation is not among the first 3
 elements. This entry is the root of the minmax tree of the permutation, and
its left subtree has at least three nodes. Let these nodes be $b_1<b_2<\cdots
<b_k$. Repeat the previous argument for this subtree, with $b_1$ playing the
role of 1 and get that whenever $b_1$ is among the leftmost three elements
of the permutation, $p_1$ is a leaf with probability 1/3. 

Iterate this algorithm as long as it is necessary. It will eventually stop 
because we either get a left subtree of size three, or a subtree whose
minimal entry is among the first three ones. In every subcase, $p_1$ is a leaf
with probability 1/3, completing the proof of the lemma. \end{proof}

We will use the same methods to prove the general theorem.
\begin{theorem} The number of $n$-permutations containing $p_i$ as a leaf 
is $n!/3$, when $1\leq i \leq n-2$. If $i=n-1$, then $p_i$ is never
a leaf, whereas if $i=n$, then $p_i$ is always a leaf.\end{theorem}

\begin{proof}
Clearly, $p_n$ is always a leaf because it cannot be the leftmost in any
 comparison,  thus it cannot have descendants. Similarly, $p_n$ is always the
 child of $p_{n-1}$, thus $p_{n-1}$ is never a leaf.

Now let $1\leq i \leq n-2$. Again, we can suppose that 1 precedes $n$ in our
 permutations.
First we consider the case when $1\in \{ p_i,p_{i+1}, p_{i+2} \} $. 
\begin{itemize}
\item If $p_i=1$, then $p_i$ is the root of the minmax tree.
\item If $p_{i+1}=1$, then $p_{i+1}$ is the root, so $p_i$ is the rightmost
element of its left subtree, and as such, it is necessarily a leaf.
\item If $p_{i+2}=1$, then $p_{i+2}$ is the root, $p_i$ is the next-to-last
element of the root's left subtree, and as such, it is always an internal node
(having the leaf $p_{i+1}$ for its only child).
\end{itemize} 
Again each of these subcases occurs with probability 1/3. 

If $1\notin \{p_i,p_{i+1},p_{i+2}\}$, then we can proceed as we did in proof of
Lemma 1: look for the entry 1 of the permutation, then only consider the
subtree which contains the positions $i, i+1$ and $i+2$. If 1 was not in any
of these positions, then all the three of them are in the same subtree. 
Iterating this algorithm we eventually reach a subtree where we can apply the
above method. The structure of the other subtrees do not influence whether
$p_i$ is a leaf or not, so $p_i$ is a leaf with probability 1/3.
\end{proof}

\begin{corollary} Let $n\geq 4$ and let $d_{i,j}$ be the number of 
$n$-permutations whose entry
$p_i$ has exactly $j$ children. Then $d_{i,j}=n!/3$ for $2\leq i\leq n-2$ and
$j\in \{0,1,2\}$. \end{corollary}

\begin{proof} For $j=0$, this is just the statement of the Theorem. For $j=2$,
note that a node $p_{i+1}$ has two children if and only if it is the local 
extremum $e_i$, and this happens exactly in the second subcase of the proof of
the theorem, thus its probability is 1/3. Similarly, for $j=1$, $p_{i+1}$
 must be
the last element of the left subtree of the local extremum $e_i$, thus $e_i=
p_{i+2}$ must hold, which is just the third subcase. \end{proof}

As $p_1$ cannot have two children, and $p_{n-1}$ cannot be a leaf, we also
have $d_{1,1}=2n!/3 $ and $d_{n-1,1}=2n!/3$. This latter is true since 
$p_{n-1}$ has
two children if and only if $p_{n-2}$ is a leaf.

We note that there exists a natural bijection between the set of minmax trees
and that of min1-min2 trees (in these trees the root is the leftmost of
the minimal and the second minimal element) and thus our results hold for
the min1-min2 trees as well.

\section*{Acknowledgement}
I am grateful to G\'abor Hetyei who introduced me into this subject.

\end{document}